\documentclass{article}
\usepackage{maa-monthly}
\usepackage{hyperref}

\theoremstyle{theorem}
\newtheorem{theorem}{Theorem}

\theoremstyle{definition}

\begin{document}

\pagebreak

\title{A new proof of the irrationality of $e$} 
\markright{Irrationality of Euler’s number}
\author{Achraf Ben Said}  

\maketitle

\begin{abstract}
In this short note, we provide a simple proof of the irrationality of Euler’s number \(e\), without using the classical estimates based on continued fractions or Taylor series expansions.
\end{abstract}

\vspace{10px}

\noindent{\textbf{Keywords:}{ Irrationality, Euler's number} 

\vspace{10px}

The purpose of this note is to present an especially short and direct proof that \( e \) is irrational. This result was first established by Euler (see \cite{Euler1}) using continued fractions. He proved that
\begin{equation*}
e=[1; 1, 2, 1, 1, 4, 1, 1, 6, \dots].
\end{equation*}
Later, simpler derivations of this continued fraction were obtained (see, for example, \cite{Cohn, Old}).

\medskip

A standard proof of the irrationality of \( e \), due to Fourier \cite{Fourier}, proceeds by contradiction and is based on the series representation
\[
e=\sum_{n=0}^\infty \frac{1}{n!}.
\]

\medskip

In 1873, Hermite significantly strengthened this result by proving the stronger statement that $e$ is transcendental (see \cite{Hermite}).

\medskip

For our purposes, we consider the sequence \( \{c(n)\}_{n \geq 0} \) (see A000522 in \cite{X}), defined recursively by
\begin{equation}
\label{recurrence}
c(0)=1, \qquad c(n)=n\,c(n-1)+1, \quad \text{for } n \geq 1.
\end{equation}

The sequence also admits the integral representation
\begin{equation}
\label{sequence}
c(n)=e\left(n!-\int_0^1 x^n e^{-x}\,dx\right),
\end{equation}
valid for every integer \( n\geq 0 \). Indeed, a direct computation shows that \( c(0)=1 \). For \( n\geq 1 \), integration by parts gives
\begin{align*}
c(n)
&= e\left(n! + e^{-1} - n\int_0^1 x^{n-1}e^{-x}\,dx\right) \\
&= 1 + ne\left((n-1)! - \int_0^1 x^{n-1}e^{-x}\,dx\right) \\
&= 1 + n\,c(n-1),
\end{align*}
which coincides with the recurrence relation in \eqref{recurrence}.

\medskip

Using \eqref{recurrence} and a straightforward induction argument, we conclude that
\begin{equation}
\label{Identity2}
c(n) \in \mathbb{N}
\end{equation}
for all natural numbers \( n \geq 0 \).

\medskip

\begin{theorem}
The number $e$ is irrational.
\begin{proof}
Let $n \geq 1$ be a natural number. Suppose, for the sake of contradiction, that $e = \frac{a}{b}$ for some $a,b \in \mathbb{N}$. Using the identity \eqref{sequence} together with \eqref{Identity2}, it follows that
\begin{equation}
\label{Identity3}
a \int_0^1 x^n e^{-x}\,dx = n!a - b\,c(n) \in \mathbb{N},
\end{equation}
for all $n \geq 0$.
On the other hand, since $e^{-x} \leq 1$ for all $x \in [0,1]$, we have
\begin{equation*}
0 < a \int_0^1 x^n e^{-x}\,dx \leq a \int_0^1 x^n \, dx = \frac{a}{n+1}.
\end{equation*}
for all $n \geq 0$. In particular, for all $n \geq a$, it holds that
\[
0 < a \int_0^1 x^n e^{-x}\,dx < 1,
\]
which contradicts \eqref{Identity3}. Therefore, $e$ is irrational.
\end{proof}
\end{theorem}

\begin{acknowledgment}{ACKNOWLEDGMENTS.}
The author thanks his high school professor Javier Rojo, as well as his PhD supervisors Santiago Boza and Javier Soria, for their guidance and support throughout his studies. He also thanks the anonymous referees for their valuable comments and clarifications, which have greatly improved the final presentation of this work. The author was partially supported by grants PID2020-113048GB-I00 and PID2024-155917NB-I00, funded by MCIN/AEI/ 10.13039/501100011033.
\end{acknowledgment}

\begin{acknowledgment}{DISCLOSURE STATEMENT.}
No potential conﬂict of interest was reported by the author.
\end{acknowledgment}

\bibliographystyle{vancouver}
\bibliography{VancouverExamples.bib}

\begin{biog}
\item[Achraf Ben Said] 
\begin{affil}
Department of Analysis and Applied Mathematics, Complutense University of
Madrid, Madrid 28040\\
achbensa@ucm.es
\end{affil}
\end{biog}

\vfill\eject

\end{document}